\documentclass[12pt]{extarticle}      
\usepackage{nicematrix,fourier,extsizes,blkarray}
\usepackage{array}
\usepackage{tikz-cd}
\newcolumntype{C}{>{$}c<{$}}

\usepackage{blkarray}

\usepackage[symbol]{footmisc}

\usepackage{amssymb}
\usepackage{eucal}
\usepackage{amsmath}
\usepackage{amsthm}
\usepackage{graphicx}
\usepackage{float}
\usepackage{framed}
\usepackage{enumerate}
\usepackage{hyperref}
\usepackage[toc,page]{appendix}

\newcommand{\R}{{\mathbb  R}}   
\newcounter{mastercounter}
\numberwithin{mastercounter}{section}

\newtheorem{thm}[mastercounter]{Theorem}

\makeatletter
\let\c@equation\c@mastercounter
\makeatother

\numberwithin{equation}{section}

\usepackage[T1]{fontenc}

\DeclareMathOperator{\tr}{tr}

\begin{document}

\title{An explicit formula for the Laplace--Beltrami operator on the Stiefel manifold
}

\author{Petre Birtea, Ioan Ca\c su, Dan Com\u{a}nescu\\
{\small Department of Mathematics, West University of Timi\c soara} 
\\
{\small Bd. V. P\^ arvan, No 4, 300223 Timi\c soara, Rom\^ania}\\
{\small Email: petre.birtea@e-uvt.ro, ioan.casu@e-uvt.ro, dan.comanescu@e-uvt.ro}}
\date{}

\maketitle

\begin{abstract}
We derive an explicit formula for the Laplace–Beltrami operator on the orthogonal Stiefel manifold, viewed as a constraint submanifold of the Euclidean space of real matrices equipped with the Frobenius metric. Using the general framework of Laplace operators on constraint manifolds, we provide the formula for the La\-place–Bel\-tra\-mi operator in terms of the ambient Euclidean coordinates. The result extends previously known cases, recovering the formulas for the sphere and the special orthogonal group as particular instances.
\end{abstract}

\maketitle

\section{Introduction}

In the last decades, there was a great interest in the study of the Laplace operator on Riemannian manifolds (frequently named Laplace-Beltrami operator), see \cite{shubin}, \cite{berger}, \cite{bakry}, \cite{svirkin}.

Let $(M,{\bf g})$ be a Riemannian manifold  (ambient space) of dimension $m$, $F=(F_1,\dots,F_k):M\to \R^k$ be a smooth set of constraint functions, with $k<m$, and $S_{\bf c}:=F^{-1}({\bf c})$, where ${\bf c}$ is a regular value for $F$. The induced Riemannian metric on $S_{\bf c}$ is  ${\bf g}_{\bf c}={\bf g}\left|_{_{TS_{\bf c}\times TS_{\bf c}}}\right.$. 

We choose $\{{\bf t}_1,\dots,{\bf t}_{m-k}\}$ an adapted local frame on $S_{\bf c}$, i.e. ${\bf t}_i\in \mathcal{X}(M)$ such that ${\bf t}_i({\bf u})\in T_{\bf u}S_{\bf c}$, for all ${\bf u}$ in the domain of a system of local coordinates $(u_1,\dots,u_m)$. Regarding ${\bf t}_i({\bf u})$ as column vectors expressed in the local frame $\frac{\partial }{\partial u_1},\dots, \frac{\partial }{\partial u_m}$, we obtain the transformation matrix 
\begin{equation}\label{transformation-matrix-11}
T:=[{\bf t}_1, \dots, {\bf t}_{m-k}]\in \mathcal{M}_{m\times (m-k)}(\R). 
\end{equation}

For $\widetilde{f}:S_{\bf c}\to \R$ a smooth function, the Laplace-Beltrami operator is given by
\begin{equation*}\label{Laplacian-Sc}
\Delta_{S_{\bf c}}\widetilde{f}({\bf u})=\tr\left(\left[{\bf g}_{\bf c}({\bf u})]^{-1}[\text{Hess}_{S_{\bf c}}\widetilde{f}\right]({\bf u})\right),\,\,\forall{\bf u}\in S_{\bf c}.
\end{equation*}
The above formula can be written using the Riemannian geometry of the ambient manifold $M$. 
Let $f:M\to \R$ be a smooth prolongation of $\widetilde{f}$, i.e. $\widetilde{f}=f_{|S_{\bf c}}$.
In \cite{laplacian}, the following formula has been proved.

\begin{thm}\label{teorema-generala}\emph{(\cite{laplacian})}
Choosing an adapted local frame  $\{{\bf t}_1,\dots,{\bf t}_{m-k}\}$, the Laplace-Beltrami operator on a constraint manifold $S_c$, written in the coordinates of the ambient space $M$, has the formula
\begin{equation*}
\Delta_{S_{\bf c}}\widetilde{f}=\tr\left((T^t[{\bf g}]T)^{-1}T^t[\emph{Hess}{f}]T\right) - \sum_{\alpha=1}^k\sigma_{\alpha}\tr\left((T^t[{\bf g}]T)^{-1}T^t\left[\emph{Hess}F_{\alpha}\right]T\right).
\end{equation*}

For the particular case when the ambient manifold is the Euclidean space, $[{\bf g}]=\mathbb{I}_m$, we have the formula
\begin{equation*}\label{laplace-general-101}
\Delta_{S_{\bf c}}\widetilde{f}=\tr\left(\left(T(T^tT)^{-1}T^t\right)[\emph{Hess}{f}]\right) - \sum_{\alpha=1}^k\sigma_{\alpha}\tr\left(\left(T(T^tT)^{-1}T^t\right)\left[\emph{Hess}F_{\alpha}\right]\right).
\end{equation*}
\end{thm}

The functions $\sigma_{\alpha}$ are the so called Lagrange multiplier functions, defined in \cite{Birtea-Comanescu-Hessian}, and are given by 
\begin{equation*}\label{sigma-101}
\sigma_{\alpha}({\bf u}):=\frac{\det \left(\text{Gram}_{(F_1,\ldots ,F_{\alpha-1},f, F_{\alpha+1},\dots,F_k)}^{(F_1,\ldots , F_{\alpha-1},F_\alpha, F_{\alpha+1},...,F_k)}({\bf u})\right)}{\det\left(\text{Gram}_{(F_1,\ldots ,F_k)}^{(F_1,\ldots ,F_k)}({\bf u})\right)},
\end{equation*}
where $\text{Gram}_{(g_1,...,g_s)}^{(f_1,...,f_r)}$ is the Gramian of the respective functions.

We apply the above general setting to the particular case of the orthogonal Stiefel manifold $St_p^n$, regarded as a submanifold  of $\mathcal{M}_{n\times p}(\R)$ endowed with the Frobenius metric.
For $n\geq p\geq 1$, we consider the orthogonal Stiefel manifold:
$$St_p^n=\{U\in \mathcal{M}_{n\times p}(\R) \,|\,U^TU=\mathbb{I}_p\}.$$
Denote with  ${\bf u}_1,...,{\bf u}_p\in \R^n$ the  vectors that form the columns of the matrix $U\in  \mathcal{M}_{n\times p}(\R)$. The condition that the matrix $U$ belongs to $St_p^n$ is equivalent with the vectors ${\bf u}_1,...,{\bf u}_p\in \R^n$ being orthonormal.

The functions that describe the constraints defining the orthogonal Stiefel manifold as a preimage of a regular value are $F_{aa},F_{bc}:\mathcal{M}_{n\times p}({\R})\rightarrow \R$ given by:
%\footnote{We denote by $\left<\cdot,\cdot\right>$ the canonical scalar product on $\R^{n}$ and by $\|\cdot \|$ the induced norm.}
\begin{equation*}
F_{aa} (U) =\frac{1}{2}\|{\bf u}_a\|^2,\,\,1 \leq a\leq p;~~~
F_{bc} (U) = \left<{\bf u}_b,{\bf u}_c\right>,\,\,1\leq b<c\leq p.  
\end{equation*}
We have ${\bf F}:\mathcal{M}_{n\times p}({\R})\rightarrow \R^{\frac{p(p+1)}{2}}$, where ${\bf F}:=\left( \dots , F_{aa},\dots ,F_{bc}, \dots\right)$ and $$St_p^n= {\bf F}^{-1}\left( \dots , \frac{1}{2},\dots ,0, \dots\right) \subset \mathcal{M}_{n\times p}({\R}).$$

In what follows, we denote by $\hbox{vec}(U)\subset \R^{np}$ the column vectorization of the matrix $U\in \mathcal{M}_{n\times p}(\R)$. 
Also, in order to not overload the notations, we make the following {convention}.
When we use the notations $\nabla f$ and $\text{Hess}f$ (with respect to the Euclidean metric on $\mathcal{M}_{n\times p}(\mathbb{R})$), we refer to them as being in {\it matrix} form, i.e., for a smooth function $f:\mathcal{M}_{n\times p}(\R)\to \R$, we denote $\widehat{f}:\R^{np}\to \R$, $\widehat{f}:=f\circ\text{vec}^{-1}$ and
$$\nabla f(U):=\text{vec}^{-1}\left(\nabla\widehat{f}(\text{vec}(U))\right),$$
$$\text{Hess}f(U):=\left[\text{Hess}\widehat{f}\right](\text{vec}(U))=\left[\begin{array}{c|c|c}
\frac{\partial^2 \widehat{f}(\text{vec}(U))}{\partial {\bf u}_1\partial {\bf u}_1} &\dots&\frac{\partial^2 \widehat{f}(\text{vec}(U))}{\partial {\bf u}_1\partial {\bf u}_p}\\
\hline
\vdots&\ddots&\vdots\\
\hline
\rule{0pt}{12pt}\frac{\partial^2 \widehat{f}(\text{vec}(U))}{\partial {\bf u}_p\partial {\bf u}_1}&\dots&\frac{\partial^2 \widehat{f}(\text{vec}(U))}{\partial {\bf u}_p\partial {\bf u}_p}\end{array}\right],$$
where
$$\frac{\partial^2 \widehat{f}}{\partial {\bf u}_i\partial {\bf u}_j}:=\left[\begin{array}{ccc}
\frac{\partial^2 \widehat{f}}{\partial {u}_{1i}\partial {u}_{1j}} &\dots&\frac{\partial^2 \widehat{f}}{\partial {u}_{1i}\partial {u}_{nj}}\\
\vdots&\ddots&\vdots\\
\rule{0pt}{12pt}\frac{\partial^2\widehat{f}}{\partial {u}_{ni}\partial {u}_{1j}}&\dots&\frac{\partial^2 \widehat{f}}{\partial {u}_{ni}\partial {u}_{nj}}\end{array}\right].$$
Also, $\Delta f(U):=\Delta\widehat{f}(\text{vec}(U))$.

For $U\in St_p^n$, we introduce the $np\times np$ matrix 
$$\Lambda(U):=\left[\begin{array}{c|c|c}
\textbf{u}_1 \textbf{u}_1^t&\dots&\textbf{u}_p \textbf{u}_1^t\\
\hline
\vdots&\ddots&\vdots\\
\hline
\rule{0pt}{12pt}\textbf{u}_1 \textbf{u}_p^t&\dots&\textbf{u}_p \textbf{u}_p^t\end{array}\right].$$

The main result of this paper, which is proved in Section 3, is the following theorem.

\begin{thm}\label{main-theorem}
The Laplace-Beltrami operator applied to a smooth function $\widetilde{f}:St_p^n\to \R$ has the formula
$${\Delta_{St^n_p} \tilde{f}(U)=\Delta f(U)-\left(n-\frac{p+1}{2}\right) \operatorname{tr}\left(U^t \nabla f(U)\right)-\frac{1}{2} \operatorname{tr}\left(\left(\mathbb{I}_p \otimes\left(U U^t\right)+\Lambda(U)\right) \emph {Hess} f(U)\right)},$$
where $f:\mathcal{M}_{n\times p}(\R)\to \R$ is a smooth prolongation of $\widetilde{f}$.
\end{thm}

An alternative formula can be obtained noticing that the matrix $\Lambda(U)$ can be written as $\Lambda(U)=K_{p,n}\cdot (U \otimes U^t)$, where $K_{p,n}$ is the commutation matrix of dimension $np\times np$ (see Lemma A.2 in \cite{eigenvalues}).

In the particular cases $p=1$ (sphere) and $p=n$ (special orthogonal group), we recover the formulas from \cite{laplacian}.

For a two-parameter family of Riemannian metrics on $St_p^n$, using the projection method described in \cite{du-2}, a formula for the Laplace-Beltrami operator has been given in \cite{du-1}. For a specific choice of these parameters, this formula is equivalent with our formula from the above theorem.

\section{The ambient geometry of the Stiefel manifold}

Using the description of the Stiefel manifold as a preimage of a regular value, 
in \cite{edelman} it is given the following elegant explicit form for the tangent space at a point $U\in St_p^n$:
$$T_U St_p^n=\{UA+(\mathbb{I}_n-UU^t)C\,|\,A\in \mathcal{M}_{p\times p}(\R),\,A=-A^t,\,C\in \mathcal{M}_{n\times p}(\R)\}.$$
%On the tangent space we consider the canonical Frobenius scalar product.
%$$\left<\Delta_1(U),\Delta_2(U)\right>=\tr(\Delta_1^T(U)\Delta_2(U)),\,\,\,\Delta_1(U),\Delta_2(U)\in T_U St_p^n.$$

In \cite{second-order}, an explicit basis for the tangent space $T_USt_p^n$ is extensively discussed and its properties are listed.
The explicit description of the vectors in this basis depends on the choice of a full rank $p\times p$ submatrix of $U$. 
Without loss of generality, we suppose that the full rank of the matrix $U$ is given by its first $p$ rows.

We split the basis, which we denote by $\mathcal{B}_U$, as the following set union of vectors
$$\mathcal{B}_U=\mathcal{B}_U^{'}\cup \mathcal{B}_U^{''}.$$

The set $\mathcal{B}_U^{'}$ is formed with tangent vectors of the form 
\begin{equation*}\label{delta-prim}
\Delta'_{ab}(U)=UA_{ab},\,\,\,1\leq a<b\leq p,
\end{equation*}
where 
$A_{ab}=(-1)^{a+b}\left({\bf f}_a {\bf f}^t_b-{\bf f}_b {\bf f}_a^t\right),\,\,\,1\leq a<b\leq p,$
form the standard basis for the $p\times p$ skew-symmetric matrices. The vectors ${\bf f}_1, \dots, {\bf f}_p$ form the canonical basis in the Euclidean space $\R^p$.

Under the above assumption about the full rank of the matrix $U$, define the tangent vectors of $\mathcal{B}_U^{''}$ as
\begin{equation*}\label{delta-secund}
\Delta''_{ic}(U)=\left(\mathbb{I}_n-UU^t\right)C_{ic},\,\,i\in \{p+1,\dots, n\},\,c\in \{1,...,p\},
\end{equation*}
where 
$C_{ic}={\bf e}_i {\bf f}^t_c$. The vectors ${\bf e}_1, \dots, {\bf e}_n$ form the canonical basis in the Euclidean space $\R^n$.\\ %\footnote{  The $n\times p$ matrix ${\bf e}_i {\bf f}^t_c$ has $1$ on the $i$-th row and $c$-th column and for the rest $0$. }$
For a fixed $c\in \{1,...,p\}$, we define:
\begin{equation*}\label{secund-3}
 _c\mathcal{B}_{U}^{''}=\left\{\Delta''_{ic}(U)\,|\,i\in \{p+1,\dots,n\}\right\}.
\end{equation*}
The set $\mathcal{B}_U^{''}$ is the set union 
%\begin{equation*}\label{secund-2}
$\mathcal{B}_U^{''}=\bigcup\limits_{c=1}^{p}{_c\mathcal{B}_{U}^{''}}.$
%\end{equation*}

We denote $Z:=\mathbb{I}_n-UU^t$ and the elements of this $n\times n$ symmetric idempotent (projection) matrix with $Z=[z_{ij}]_{1\leq i,j\leq n}$.

The following result summarizes the properties of the basis $\mathcal{B}_U$.

\begin{thm}\label{properties-basis} \emph{(\cite{second-order})}
For $U\in St^n_p$, the set $\mathcal{B}_U=\mathcal{B}_U^{'}\cup \mathcal{B}_U^{''}$
%, where the elements of the set $\mathcal{B}_U^{'}$ are given by \eqref{delta-prim} and the elements of the set $\mathcal{B}_U^{''}$ are given by 
%\eqref{delta-secund}, \eqref{secund-3}, and \eqref{secund-2},  
is a basis for the tangent space $T_U St_p^n$. The vectors of this basis have the following properties:
\begin{itemize}
\item [(i)] $\mathcal{B}_U^{'}$ is an orthogonal set;
\item [(ii)] For $c_1,c_2\in \{1,...,p\}$ and $c_1\neq c_2$, we have $_{c_1}\mathcal{B}_{U}^{''}\perp {_{c_2}\mathcal{B}_{U}^{''}}$;
\item [(iii)] $\mathcal{B}_U^{'}\perp \mathcal{B}_U^{''}$;
\item[(iv)] For $1\leq a<b\leq p$, we have $\left<\Delta'_{ab}(U),\Delta'_{ab}(U)\right>=2$;
\item[(v)] For $k_1,k_2\in \{p+1,\dots,n\},~ c\in \{1,\dots, p\}$, we have $\left<\Delta''_{k_1c}(U),\Delta''_{k_2c}(U)\right>=z_{k_1k_2}$.
\end{itemize}
\end{thm}

\section{Proof of the main result}

For the matrices $U\in St_p^n$ and $Z=\mathbb{I}_n-UU^t$, we consider the two-block column decompositions
\begin{equation*}
U=\left[\begin{array}{c}U_1 \\ \hline U_2\end{array}\right],~ U_1 \in \mathcal{M}_{p \times p},~ U_2 \in \mathcal{M}_{(n-p) \times p};~~~Z=\left[\begin{array}{c}Z_1 \\ \hline Z_2\end{array}\right],~Z_1 \in \mathcal{M}_{p \times n},~ Z_2 \in \mathcal{M}_{(n-p) \times n}.
\end{equation*}
From the definition of a Stiefel matrix $\left(U^t U=\mathbb{I}_p\right)$, it follows that 
\begin{equation}\label{orto}
U_1^t U_1+U_2^t U_2=\mathbb{I}_p.
\end{equation}
We also have
$$Z=\left[\begin{array}{cc}
    \mathbb{I}_p-U_1U_1^t & -U_1U_2^t \\
    -U_2U_1^t & \mathbb{I}_{n-p}-U_2U_2^t
\end{array}
\right],$$
hence
\begin{equation}\label{Z2}
Z_2=\left[\begin{array}{ccc} -U_2U_1^t&|& \mathbb{I}_{n-p}-U_2U_2^t \end{array}\right].
\end{equation}
From the equality $ZZ^t=Z$ it follows that
\begin{equation}\label{Z2-Z2}
    Z_2Z_2^t=\mathbb{I}_{n-p}-U_2U_2^t.
\end{equation}
We organize the transformation matrix $T\in \mathcal{M}_{np\times\left( \frac{p(p-1)}{2}+p(n-p)\right)}(\R)$  in the point $U$ (see \eqref{transformation-matrix-11}) as
$$T=\left[{\bf t}'_{12} ,\dots, {\bf t}'_{1p},{\bf t}'_{23},\dots, {\bf t}'_{2p},\dots, {\bf t}'_{p-1,p},{\bf t}''_{p+1,1},\dots,{\bf t}''_{n,1},\dots, {\bf t}''_{p+1,p},\dots, {\bf t}''_{n,p}\right],$$
where
${\bf t}'_{ab}(U)=\text{vec}\left(\Delta'_{ab}(U)\right)$ and ${\bf t}''_{ic}(U)=\text{vec}\left(\Delta''_{ic}(U)\right)$ for $1\leq a<b\leq p$ and $i\in \{p+1,\dots, n\},\,c\in \{1,...,p\}$. The vectorization of the elements from $\mathcal{B}'_U$ are in lexicographic order of indices, while that of the elements from $\mathcal{B}''_U$ are in colexicographic order of indices.

If we decompose the matrix $T$ in a two-block row (corresponding to the vectors in the two parts of the basis for the tangent space, $\mathcal{B}_U'$ and $\mathcal{B}_U''$ respectively), i.e.
\begin{equation}\label{descompunere-T}
T=\left[~ T_1 ~| ~T_2~\right],
\end{equation}
with $T_1\in \mathcal{M}_{np \times \frac{p(p-1)}{2}}$ and  $T_2\in \mathcal{M}_{np \times p(n-p)}$, then
\begin{equation}\label{Ttranspus-T}
T^tT=\left[
\begin{array}{cc}
    T_1^tT_1 & T_1^tT_2 \\
    ~\\
    T_2^tT_1 & T_2^tT_2
\end{array}
\right].
\end{equation}
From Theorem \ref{properties-basis} (i) and (iv) it follows that 
$T_1^tT_1=2 \mathbb{I}_{\frac{p(p-1)}{2}}.$
From Theorem \ref{properties-basis} (iii) it follows that
$T_1^tT_2=\mathbb{O}_{\frac{p(p-1)}{2}\times p(n-p)}$
and
$T_2^tT_1=\mathbb{O}_{p(n-p)\times \frac{p(p-1)}{2}}.$\\
\vspace*{0.1cm}\\
From the definition of the set $\mathcal{B}''_U$ and of the vectors ${\bf t}''_{ic}(U)$ it follows that\footnote{$\text{vec}(AXB)=(B^t\otimes A)\text{vec}(X)$ (see \cite{cookbook}).}
$${\bf t}''_{ic}=\hbox{vec}\left(Z{\bf e}_i {\bf f}^t_c\right)=\left({\bf f}_c\otimes Z\right){\bf }\hbox{vec}({\bf e}_i)=\left({\bf f}_c\otimes Z\right){\bf }{\bf e}_i=\left[\begin{array}{c}
\mathbb{O}_{n\times 1}\\
\vdots\\
 Z{\bf e}_i\\
\vdots\\
\mathbb{O}_{n\times 1}
\end{array}\right]=\begin{blockarray}{c c}
  % & C_1 & C_2 & C_3 \\
\begin{block}{[c] c}
\mathbb{O}_{n\times 1} &  \\
\vdots & \\
 {\bf z}_i & \leftarrow c^{\hbox{th}} \hbox{ line}\\
\vdots & \\
\mathbb{O}_{n\times 1}&  \\
\end{block}
\end{blockarray}.$$
Therefore, it is easy to notice that
$$T_2=\mathbb{I}_p\otimes Z_2^t,$$
hence
$T_2^t=\mathbb{I}_p\otimes Z_2$
and, see \eqref{Z2-Z2},
$$T_2^tT_2=\mathbb{I}_p\otimes \left(Z_2Z_2^t\right)=\mathbb{I}_p\otimes \left(\mathbb{I}_{n-p}-U_2U_2^t\right).$$
Combining the above equalities in \eqref{Ttranspus-T}, we obtain
\begin{equation*}
T^t T=\left[\begin{array}{cc}
2 \mathbb{I}_{\frac{p(p-1)}{2}} & \mathbb{O}_{\frac{p(p-1)}{2}\times p(n-p)} \\
~\\
\mathbb{O}_{p(n-p)\times \frac{p(p-1)}{2}} & \mathbb{I}_p \otimes \left(\mathbb{I}_{n-p}-U_2U_2^t\right)  
\end{array}\right] 
\end{equation*}
and therefore
\begin{equation}\label{Ttranspus-T-invers}
(T^t T)^{-1}=\left[\begin{array}{cc}
\frac{1}{2} \mathbb{I}_{\frac{p(p-1)}{2}} & \mathbb{O}_{\frac{p(p-1)}{2}\times p(n-p)} \\
~\\
\mathbb{O}_{p(n-p)\times \frac{p(p-1)}{2}} & \mathbb{I}_p \otimes \left(\mathbb{I}_{n-p}-U_2U_2^t\right)^{-1}  
\end{array}\right] .
\end{equation}
By the Woodbury matrix inversion formula and \eqref{orto}, we have:
\begin{equation}\label{woodbury}
\begin{aligned}
\left(\mathbb{I}_{n-p}-U_2U_2^t\right)^{-1}&=\mathbb{I}_{n-p}+U_2\left(\mathbb{I}_p-U_2^t U_2\right)^{-1} U_2^t\\
&=\mathbb{I}_{n-p}+U_2\left(U_1^t U_1\right)^{-1} U_2^t\\
&=\mathbb{I}_{n-p}+U_2U_1^{-1}\left(U_1^t\right)^{-1} U_2^t.
\end{aligned}
\end{equation}
From \eqref{descompunere-T}, \eqref{Ttranspus-T-invers} and \eqref{woodbury}
we obtain
\begin{equation}\label{product-13}
\begin{aligned}
T(T^tT)^{-1}T^t&=
\dfrac{1}{2}T_1T_1^t+
T_2\left(\mathbb{I}_p\otimes \left(\mathbb{I}_{n-p}+U_2U_1^{-1}\left(U_1^t\right)^{-1} U_2^t\right) \right)T_2^t\\
&=\dfrac{1}{2}T_1T_1^t+\mathbb{I}_p\otimes\left(Z_2^t\left(\mathbb{I}_{n-p}+U_2U_1^{-1}\left(U_1^t\right)^{-1} U_2^t\right) Z_2\right).
\end{aligned}
\end{equation}
Closely following the proof of Lemma 3.2 from \cite{laplacian}, we have
\begin{equation}\label{T1-35}
T_1T_1^t=\mathbb{I}_p\otimes (UU^t)-\Lambda(U).
\end{equation}
By a straightforward computation, using \eqref{orto} and \eqref{Z2}, we obtain
$$
Z_2^t\left(\mathbb{I}_{n-p}+U_2U_1^{-1}\left(U_1^t\right)^{-1} U_2^t\right) Z_2= \left[\begin{array}{c}-U_1U_2^t \\  \mathbb{I}_{n-p}-U_2U_2^t\end{array}\right] \left(\mathbb{I}_{n-p}+U_2U_1^{-1}\left(U_1^t\right)^{-1} U_2^t\right) Z_2=$$
\begin{align*}
&=\left[\begin{array}{c}-U_1U_2^t-U_1(\mathbb{I}_p-U_1^tU_1)U_1^{-1}(U_1^t)^{-1}U_2^t \\  \mathbb{I}_{n-p}-U_2U_2^t+U_2U_1^{-1}\left(U_1^t\right)^{-1} U_2^t-U_2(\mathbb{I}_p-U_1^tU_1)U_1^{-1}(U_1^t)^{-1}U_2^t
\end{array}\right] Z_2\\
&=\left[\begin{array}{c} -(U_1^t)^{-1}U_2^t\\  \mathbb{I}_{n-p}
\end{array}\right]\cdot \left[\begin{array}{ccc} -U_2U_1^t&|& \mathbb{I}_{n-p}-U_2U_2^t \end{array}\right] \\
&=\left[\begin{array}{cc}
   (U_1^t)^{-1}(\mathbb{I}_p-U_1^tU_1)U_1^t  & -(U_1^t)^{-1}U_2^t+(U_1^t)^{-1}(\mathbb{I}_p-U_1^tU_1)U_2^t \\
   -U_2U_1^t  & \mathbb{I}_{n-p}-U_2U_2^t
\end{array}\right]\\
&=\left[\begin{array}{cc}
    \mathbb{I}_p-U_1U_1^t & -U_1U_2^t \\
    -U_2U_1^t & \mathbb{I}_{n-p}-U_2U_2^t
\end{array}
\right].\end{align*}
Therefore,
\begin{equation}\label{Z2-123}
  Z_2^t\left(\mathbb{I}_{n-p}+U_2U_1^{-1}\left(U_1^t\right)^{-1} U_2^t\right) Z_2=Z=\mathbb{I}_n-UU^t.  
\end{equation}
Combining now equations \eqref{product-13}, \eqref{T1-35}, \eqref{Z2-123} we obtain
\begin{equation*}
T(T^tT)^{-1}T^t=\frac{1}{2}\mathbb{I}_p\otimes \left(UU^t\right)-\frac{1}{2}\Lambda(U)+\mathbb{I}_p\otimes \left(\mathbb{I}_n-UU^t\right),
\end{equation*}
or, equivalently,
\begin{equation}\label{produs-T-121}
T(T^tT)^{-1}T^t=
\mathbb{I}_{np}-\frac{1}{2}\mathbb{I}_p\otimes (UU^t)-\frac{1}{2}\Lambda(U).
\end{equation}

The Lagrange multiplier functions for the case of the orthogonal Stiefel manifold, see \cite{first-order} and \cite{second-order}, are given by the formulas: 
\begin{equation}\label{Lagrange-multipliers-functions}
\sigma_{aa}(U)=\biggl \langle\displaystyle\frac{\partial f}{\partial {\bf u}_a}(U),{\bf u}_a\biggr \rangle;~~
\sigma_{bc}(U)=\displaystyle\frac{1}{2}\left(\biggl \langle\displaystyle\frac{\partial f}{\partial {\bf u}_c}(U),{\bf u}_b\biggr \rangle+\biggl \langle\displaystyle\frac{\partial f}{\partial {\bf u}_b}(U),{\bf u}_c\biggr \rangle\right).
\end{equation}

Considering the symmetric matrix
$\Sigma(U): =\left[\sigma_{bc}(U)\right]\in \mathcal{M}_{p\times p}(\R),$
where $\sigma_{cb}(U):=\sigma_{bc}(U)$ for $1\leq b<c\leq p$.
Using \eqref{Lagrange-multipliers-functions}, we have (see \cite{second-order})
\begin{equation}\label{Sigma-matriceal}
\Sigma(U)=\frac{1}{2}\left(\nabla f(U)^tU+U^t\nabla f(U)\right).
\end{equation}

For the case of $St_p^n$, the formula for the Hessian matrix of the function $\widetilde{f}$, as given in \cite{second-order}, is
\begin{equation*}\label{}
\hbox{Hess}_{St_p^n}\,\widetilde{f}(U) =\left(\hbox{Hess}\,f(U)-{\Sigma}(U)\otimes \mathbb{I}_n\right)_{|T_U St_p^n\times T_U St_p^n}.
\end{equation*}
We can now proceed to deduce the formula for the Laplace-Beltrami operator announced in Theorem \ref{main-theorem}. More precisely, from \eqref{produs-T-121}, we successively obtain
\begin{equation}\label{prefinal}
\begin{aligned}
\Delta_{St_p^n}\tilde{f}(U)&=\tr\left(T(T^tT)^{-1}T^t  \left(\hbox{Hess}\,f(U)-{\Sigma}(U)\otimes \mathbb{I}_n\right)\right)\\
&=\tr\left(\left(\mathbb{I}_{np}-\frac{1}{2}\mathbb{I}_p\otimes (UU^t)-\frac{1}{2}\Lambda(U)\right)\left(\hbox{Hess}\,f(U)-{\Sigma}(U)\otimes \mathbb{I}_n\right)\right)\\
&=\Delta f(U)-\frac{1}{2} \operatorname{tr}\left(\left(\mathbb{I}_p \otimes\left(U U^t\right)+\Lambda(U)\right) \text {Hess} f(U)\right)-\tr\left(\Sigma(U)\otimes \mathbb{I}_n\right)+\\
&+\frac{1}{2}\tr\left(\Sigma(U)\otimes (UU^t)\right)+
\frac{1}{2}\tr\left(\Lambda(U)\cdot \left(\Sigma(U)\otimes \mathbb{I}_n\right) \right).
\end{aligned}
\end{equation}
From \eqref{Sigma-matriceal} it follows that $$\tr\left(\Sigma(U)\otimes \mathbb{I}_n\right)=\tr(\Sigma(U))\cdot \tr(\mathbb{I}_n)=n\tr(U^t\nabla f(U))$$ 
and
$$\tr\left(\Sigma(U)\otimes (UU^t)\right)=\tr(\Sigma(U))\cdot \tr(UU^t)=\tr(\Sigma(U))\cdot \tr(U^tU)=p\tr(U^t\nabla f(U)).$$
As in the proof of Theorem 3.3 in \cite{laplacian}, we easily derive that
$$\tr\left(\Lambda(U)\cdot \left(\Sigma(U)\otimes \mathbb{I}_n\right) \right)=\tr(\Sigma(U))=\tr(U^t\nabla f(U)).$$
Combining the above equalities in \eqref{prefinal}, we obtain the formula from the Theorem \ref{main-theorem}. \qedsymbol

\end{document}